\documentclass[12pt]{article}

\usepackage{amsmath,amsthm,amsfonts,amssymb,amscd,cite}
\usepackage{tikz}
\usepackage{color}
\usepackage{hyperref}

\newtheorem{claim}{Claim}
\newtheorem{theorem}{Theorem}[section]

\newtheorem{corollary}[theorem]{Corollary}

\date{}
\begin{document}
	
\title{On Neutral Edge Sets in Anti-Ramsey Numbers}
	
\author{
Ali Ghalavand$^{a,}$\thanks{ email:\texttt{alighalavand@nankai.edu.cn }} \and  Qing Jie$^{b,}$\thanks{ email:\texttt{qingjie1126@163.com}} \and  Zemin Jin$^{b,}$\thanks{Corresponding author, email:\texttt{zeminjin@zjnu.cn }}\and  Xueliang Li$^{a,c,}$\thanks{ email:\texttt{lxl@nankai.edu.cn }}\and Linshu  Pan $^{b,}$\thanks{ email:\texttt{1925984948@qq.com }} }
	
\maketitle
	
\begin{center}
$^a$ Center for Combinatorics and LPMC, Nankai University, Tianjin 300071, China \\
$^b$ School of Mathematical Sciences, Zhejiang Normal University	Jinhua 321004, China\\
$^c$ School of Mathematical Sciences, Xinjiang Normal University, Urumchi, Xinjiang, China
\end{center}
	
\begin{abstract}
The anti-Ramsey number of a graph $G$, introduced by Erdős et al.\ in 1975, is the maximum number of colors in an edge-coloring of the complete graph $K_n$ that avoids a rainbow copy of $G$. We call a subset of edges of $G$ \emph{neutral} for the anti-Ramsey number if removing them does not alter the anti-Ramsey number of $G$. Let $k$, $t$, and $n$ be positive integers, and consider $G = kP_4 \cup tP_2$. Assume $S \subseteq E(G)$ consists of internal edges of the $P_4$ components in $G$. It is known that $S$ is neutral when $t \geq k+1 \geq 2$ and $n \geq 8k + 2t - 4$. In this paper, we identify values of $k \geq t$ such that, for all $n$ in a specific subinterval of $[8k + 2t - 4, \infty)$, $S$ remains neutral. Since the anti-Ramsey numbers for matchings are well understood, our results provide a complete determination of the anti-Ramsey number for $G$ under these conditions. Based on our findings, we conjecture that this neutrality may extend to the general case $t \geq 1$, $k \geq 1$, and $n \geq 4k + 2t$, but not when $t = 0$, $k \geq 2$, and $n \geq 4k$.
\end{abstract}

	
	
\section{Introduction}
This study investigates the effect of edge deletion on the anti-Ramsey numbers for linear forests consisting of paths of order two and four (namely, \( P_2 \) and \( P_4 \)). We begin by reviewing key definitions and providing some historical context.
Let \( G \) be a finite, simple graph. An \emph{edge-coloring} of \( G \) is a function \( c : E(G) \to \mathcal{C} \), where \( \mathcal{C} \) is a set of colors. A subgraph \( H \subseteq G \) is called \emph{rainbow} if all its edges receive distinct colors under \( c \). For a subgraph \( H \) or an edge set \( E' \subseteq E(G) \), we denote the set of colors used on its edges by \( c(H) \) or \( c(E') \), respectively.
Given a vertex subset \( V_1 \subseteq V(G) \), the \emph{induced subgraph} \( G[V_1] \) is the subgraph with vertex set \( V_1 \) and all edges of \( G \) whose endpoints are both in \( V_1 \). For another subset \( V_2 \subseteq V(G) \), the notation \( [V_1, V_2]_G \) refers to the set of edges in \( G \) with one endpoint in \( V_1 \) and the other in \( V_2 \). The \emph{disjoint union} of two graphs \( G \) and \( G' \) is denoted by \( G \cup G'\). For a positive integer \( k \), \( kG \) denotes the disjoint union of \( k \) copies of \( G \). If \( E' \subseteq E(G) \), then \( G - E' \) denotes the graph obtained by removing the edges in \( E' \) from \( G \).
For a positive integer \( \nu \), we let \( [\nu] = \{1, 2, \ldots, \nu\} \), and \( [0] = \emptyset \). The complete graph, cycle, and path on \( \iota \) vertices are denoted by \( K_\iota \), \( C_\iota \), and \( P_\iota \), respectively.
 
Given a positive integer \( n \) and a graph \( G \), the \emph{anti-Ramsey number} \( AR(n, G) \) is defined as the maximum number of colors in an edge-coloring of the complete graph \( K_n\) that contains no rainbow copy of \( G \). Introduced by Erdős et al.\ in 1975 \cite{Er1}, this concept is closely related to Turán numbers, and their seminal work included conjectures for the anti-Ramsey numbers of paths and cycles.
Subsequent research has thoroughly explored anti-Ramsey numbers for edge-disjoint subgraphs and for graphs consisting of several small disjoint components. Results concerning anti-Ramsey numbers for spanning trees in complete graphs and more general structures can be found in \cite{R2-SP-1,R2-SP-2,p1-9,R2-SP-3,R2-SP-4}. Anti-Ramsey numbers for matchings (collections of disjoint edges) in complete graphs are addressed in \cite{p1-2, p1-5, p1-7, p1-9, p1-14}, and extensions to uniform hypergraphs appear in \cite{R2-HP-1,R2-HP-2,R2-HP-3}.
More recently, attention has focused on more complex structures, particularly vertex-disjoint unions of small connected graphs and matchings. Our work contributes to this ongoing investigation by studying linear forests consisting of short odd paths. For clarity, we summarize below the principal results most relevant to our study.

Simonovits and Sós \cite{p2-10} determined $AR(n, P_t)$ for sufficiently large $t$ and $n$, while Montellano-Ballesteros and Neumann-Lara \cite{p3-6} provided a formula for $AR(n, C_t)$ when $n \geq t \geq 3$. Schiermeyer \cite{p2-9} established $AR(n, tP_2)$ for all $t \geq 2$ and $n \geq 3t + 3$. Chen et al. \cite{p1-2} and Fujita et al. \cite{p1-5} showed that for all $t \geq 2$ and $n \geq 2t + 1$, 
\begin{equation}\label{eqth0} 
AR(n, tP_2) = \begin{cases} (t-2)(2t-3) + 1 & \text{if } n \leq \frac{5t-7}{2}, \\
(t-2)\left(n - \frac{t-1}{2}\right) + 1 & \text{if } n \geq \frac{5t-7}{2}.
\end{cases} 
\end{equation} 
The remaining case $n = 2t$ was addressed by Haas and Young \cite{p1-7}, who confirmed the conjecture in \cite{p1-5}, thus completing the determination of $AR(n, tP_2)$ for all $n$ and $t$.
Bialostocki et al. \cite{p1-1} determined the anti-Ramsey numbers for $P_4$ and $P_4 \cup P_2$, and also proved that $AR(n, P_3 \cup P_2) = 2$ for $n \geq 5$, and $AR(n, P_3 \cup 2P_2) = n$ for $n \geq 7$. For broader classes, Gilboa and Roditty \cite{p1-6} established upper bounds for anti-Ramsey numbers of graphs such as $L \cup tP_2$ and $L \cup tP_3$ under certain conditions, and provided exact results for $AR(n, P_{k+1} \cup tP_3)$, $AR(n, kP_3 \cup tP_2)$, $AR(n, P_2 \cup tP_3)$, $AR(n, P_3 \cup tP_2)$, $AR(n, P_4 \cup tP_2)$, $AR(n, C_3 \cup tP_2)$, and $AR(n, tP_3)$ for sufficiently large $n$.
He and Jin \cite{p1-8} determined $AR(n, P_3 \cup tP_2)$ for $t \geq 2$ and $n \geq 2t + 3$, and also calculated $AR(n, 2P_3 \cup tP_2)$ for $t \geq 2$ and $n \geq 2t + 7$. Jie et al. \cite{p1-10} studied $AR(n, kP_3 \cup tP_2)$, providing results when $k \geq 2$, $t \geq \frac{k^2 - k + 4}{2}$, and $n \geq 3k + 2t + 1$. Recently, two of the present authors \cite{new-2} obtained a formula for $AR(n, kP_3 \cup tP_2)$ under the conditions $k \geq 2$, $t \geq \frac{k^2 - 3k + 4}{2}$, and $n = 2t + 3k$, which was further generalized by other authors in \cite{GL-1}.
Jin and Gu \cite{new-1} investigated the anti-Ramsey number for graphs whose components are $K_4$ or $P_2$, establishing $AR(n, K_4 \cup tP_2)$ for all $n \geq \max\{7, 2t + 4\}$ and $t \geq 1$. In addition, the authors of \cite{new-3} determined $AR(n, P_5 \cup tP_2)$ for linear forests of $P_2$ and $P_5$ components.
Fang et al. \cite{p1-4} considered linear forests $F$ with components of orders $p_1, p_2, \ldots, p_t$ (with $t \geq 2$ and $p_i \geq 2$ for all $i$, and at least one $p_i$ even), determining $AR(n, F)$ for sufficiently large $n$. Building upon this, Xie et al. \cite{p1-14} established an exact formula for anti-Ramsey numbers of linear forests with at least one even component and sufficiently large $n$.
Very recently, the present authors \cite{GL-2} studied the effect of edge deletion on anti-Ramsey numbers for linear forests composed of $P_2$ and $P_4$. They showed that deleting certain edges may decrease or leave unchanged the anti-Ramsey number. Specifically, let $G=kP_4 \cup tP_2$ and let $E'$ be a subset of $E(G)$ such that each endpoint of every edge in $E'$ has degree two in $G$. They proved that if any of the conditions (i) $t \geq k+1 \geq 2$ and $n \geq 8k + 2t - 4$; (ii) $k, t \geq 1$ and $n = 4k + 2t$; (iii) $k = 1$, $t \geq 1$, and $n \geq 2t + 4$ hold, then $AR(n,G) = AR(n, G-E')$. However, this does not hold when $k \geq 2$, $t = 0$, and $n = 4k$.
As a consequence, they determined $AR(n,kP_4 \cup tP_2)$ for the cases: (i) $t \geq k+1 \geq 2$ and $n \geq 8k + 2t - 4$; (ii) $k, t \geq 1$ and $n = 4k + 2t$; (iii) $k = 1$, $t \geq 0$, and $n \geq 2t + 4$; (iv) $k \geq 1$, $t = 0$, and $n = 4k$. The following theorem from this study plays a crucial role in our proof.
\begin{theorem}\label{th1}{\rm\cite{GL-2}}
For three integers $k$, $ t $, and $ n $, if  $ t\geq k+1\geq 2 $, and $ n \ge 8k+2t-4 $, then
	\[
AR(n, kP_4 \cup tP_2) = AR(n, (2k+t)P_2).
		\]
\end{theorem}

Let $G = kP_4 \cup tP_2$, and let $E'$ be the subset of the edge set of $G$ defined above. In this paper, we show that for all integers $t \geq 2$, $k \geq 1$, and $n\geq0$, if $k\geq t$ and
\[ n \in \left[8k + 2t - 4,\frac{7k^{2} + 8kt - t^{2} + 23k - t - 18}{2k - 2t + 2}\right], \] then it holds that $AR(n, G) = AR(n, G - E')$. In particular, we establish the following theorem:

\begin{theorem}\label{th2} 
Let $k$, $t$, and $n$ be integers with $t \geq 2$, $k \geq 1$, $n\geq0$. If $k\geq t$ and \[ n \in \left[8k + 2t - 4, \frac{7k^{2} + 8kt - t^{2} + 23k - t - 18}{2k - 2t + 2}\right], \] then \[ AR(n, kP_4 \cup tP_2) = AR(n, (2k+t)P_2). \]
\end{theorem}

Note that, since $k \geq t$, the interval $\left[8k + 2t - 4, \frac{7k^{2} + 8kt - t^{2} + 23k - t - 18}{2k - 2t + 2}\right]$ is nonempty if $k$ belongs to the interval $\left[t, \frac{10 t}{9} + \frac{5}{6} + \frac{\sqrt{508 t^{2} + 132 t -135}}{18}\right]$.
The following corollary follows directly from Theorem~\ref{th2} and Equation~\eqref{eqth0}
\begin{corollary}
Let $k$, $t$, and $n$ be integers with $t \geq 2$, $k \geq 1$, $n\geq0$. If $k\geq t$ and \[ n \in \left[8k + 2t - 4, \frac{7k^{2} + 8kt - t^{2} + 23k - t - 18}{2k - 2t + 2}\right], \] then
	\[AR(n,kP_4\cup tP_2) =(2k + t - 2)\left(n - \frac{2k+t-1}{2}\right) + 1.\]
\end{corollary}

\section{Proof of Theorem~\ref{th2}}
Let \( k \), \( t \), and \( n \) be positive integers with \( k \geq t \geq 2 \), and assume \[ n \in \left[8k + 2t - 4, \frac{7k^{2} + 8kt - t^{2} + 23k - t - 18}{2k - 2t + 2}\right]. \] Since every rainbow path \( P_4 \) contains a rainbow subgraph isomorphic to \( 2P_2 \), it follows that \[ AR(n, kP_4 \cup tP_2) \geq AR(n, (2k + t)P_2). \] To prove our theorem, it suffices to establish the reverse inequality: 
\begin{equation}\label{seq000} 
AR(n, kP_4 \cup tP_2) \leq AR(n, (2k + t)P_2).
\end{equation} 
To do this, we will show that if \( c \) is any edge-coloring of \( K_n \) with \( |c(K_n)| = AR(n, (2k+t)P_2) + 1 \), then \( K_n \) must contain a rainbow subgraph isomorphic to \( kP_4 \cup tP_2 \).
We proceed by first establishing that if \( K_n \) contains a rainbow subgraph isomorphic to \( (k-1)P_4 \cup (t+2)P_2 \), then it also contains a rainbow subgraph isomorphic to \( kP_4 \cup tP_2 \).

Suppose, for the sake of contradiction, that \( K_n \) does not contain a rainbow subgraph isomorphic to \( kP_4 \cup tP_2 \), but does contain a rainbow subgraph \( H \cong (k-1)P_4 \cup (t+2)P_2 \).
Let \( H_1 = (k-1)P_4 = \{ P_4^i : i \in [k-1] \} \) and \( H_2 = (t+2)P_2 = \{ P_2^j : j \in [t+2] \} \) denote the two components of \( H \). The edges of each \( P_4^i \) are given by \( E(P_4^i) = \{ x_1^i x_2^i, x_2^i x_3^i, x_3^i x_4^i \} \) for \( i \in [k-1] \), and each \( P_2^j \) has edge \( E(P_2^j) = \{ y_1^j y_2^j \} \) for \( j \in [t+2]\).
Let \( G \) be a rainbow spanning subgraph of \( K_n \) with \( |c(K_n)| \) edges, containing all edges of \( H \) (i.e., \( E(H) \subseteq E(G) )\). Based on these assumptions, we can make several structural claims about \( G \).
Before presenting the first claim, we introduce two sets, $\mathcal{X}(G)$ and $\mathcal{Y}(G)$, which together partition the edge set of $G$ as follows:
\begin{align*} 
 \mathcal{X}(G) &= E(G[V(H_2)]) \cup E(G[V']) \cup [V(H_2), V']_G, \\
  \mathcal{Y}(G) &= E(G[V(H_1)]) \cup [V(H_1), V']_G \cup [V(H_1), V(H_2)]_G.
\end{align*}

\begin{claim}\label{scl0} 
It holds that $| \mathcal{X}(G)| \leq 2n -8k- 3t +2$. 
\end{claim}

To prove this claim, By our assumptions, for each $i \neq j \in [t+2]$, we have $[V(P_2^i), V(P_2^j)]_G = \emptyset$. Thus, $|E(G[V(H_2)])|=t+2$. Furthermore, any component of $G[V']$ that is not a star must be isomorphic to $K_1$, $K_2$, or $K_3$. For indices $i \neq j \in [t+2]$ and vertices $u \neq v \in V'$, the subgraph $G[{y^i_1, y^i_2, y^j_1, y^j_2, u}]$ always contains at least one isolated edge.
Additionally, for any $i \in [t+2]$, $j \in [2]$, and $v \in V'$, if $y^i_j v \in E(G)$, then there does not exist $u \in (V' - {v})$ such that $uv \in E(G)$.
Therefore, the total number of edges in $E(G[V']) \cup [V(H_2), V']_G$ is at most $2|V'| = 2n-8k - 4t$.
Combining this with the bound for $|E(G[V(H_2)])|$, we obtain the desired inequality: $| \mathcal{X}(G)| \leq 2n-8k - 3t +2$.   This completes the proof. \qed
	
Our next objective is to establish an upper bound for the size of the set $\mathcal{Y}(G)$. To do this, we first prove several intermediate claims. Before stating the next claim, we introduce some additional definitions and notation. Let $\Gamma_1$ denote the set of all paths of length four that comprise the graph $H_1$. For each \( h \in [k-1] \), we define the subset \(\mathcal{F}_h \subseteq \Gamma_1\) as follows: \\
{\bf 1.} The set \(\mathcal{F}_1\) is the largest subset of \(\Gamma_1\) such that for any two elements \(P^i_4\) and \(P^j_4\) in this set, the condition \(|[V(P^i_4),V(P^j_4)]_G| \geq 14\) holds.\\
{\bf 2.} For \(h\) ranging from \(2\) to \(k-1\), define \(\Gamma_h = \Gamma_{h-1} - \mathcal{F}_{h-1}\) and consider \(\mathcal{F}_h\) as the largest subset of \(\Gamma_h\) such that for any two elements \(P^i_4\) and \(P^j_4\) in this set, the condition \(|[V(P^i_4),V(P^j_4)]_G| \geq 14\) holds.
	
Additionally, let \(\lambda\) be the largest integer in \([k-1]\) such that \(\mathcal{F}_\lambda \neq \emptyset\). For \(h \in [\lambda]\), define \(\mathcal{V}_h = \bigcup_{P_4 \in \mathcal{F}_h} V(P_4)\). Based on our definitions, for any two distinct paths \(P^i_4\) and \(P^j_4\) in the set \(\mathcal{F}_h\), we have 
\begin{equation}\label{seq1}
	|[V(P^i_4),V(P^j_4)]_G| \geq 14.
\end{equation} 
	
\begin{claim}\label{scl1}
Let \(h \in [\lambda]\), $i\in[t+2]$, $i_1\in[2]$, and $v_1\neq v_2\in V'$. If $|\mathcal{F}_h|\geq2$, then there are not three distinct vertices $u_j$, $j\in [3]$, in $\mathcal{V}_h$ such that $\{v_1u_1,v_2u_2,y^i_{i_1}u_3\}\subseteq E(G)$. 
\end{claim}

To prove this claim, we will assume, for the sake of contradiction, that \( |\mathcal{F}_h| \geq 2 \) and that there exist three distinct vertices \( u_j \), where \( j \in [3] \), in \( \mathcal{V}_h \) such that \( \{v_1u_1, v_2u_2, y^i_{i_1}u_3\} \subseteq E(G) \). We will consider three non-equivalent cases based on the occurrences of the vertices \( u_j \).\\
{\bf Case 1.} For an element \( P^l_4 \) in the set \( \mathcal{F}_h \), assume that \( \{u_i : i \in [3]\} \subseteq V(P^l_4) \). Let us define \( V(P^l_4) - \{u_i : i \in [3]\} = \{x^l_{l_1}\} \), and let \( P^r_4 \) be a different element from the set \( \mathcal{F}_h \). Under these circumstances, we can derive the following propositions from Equation \eqref{seq1}:\\
\noindent
{\bf 1.1} If \( u_1u_2 \) and \( x^l_{l_1}u_3 \in E(P^l_4) \), then the paths \( v_1u_1u_2v_2 \) and \( y^i_{i_2}y^i_{i_1}u_3x^l_{l_1} \), where \( i_2 \in ([2] - \{i_1\}) \), create a rainbow subgraph of \( G \) that is isomorphic to \( 2P_4 \).\\
\noindent
{\bf 1.2} If \( u_1u_2 \in E(P^l_4) \) and \( x^l_{l_1}u_3 \notin E(P^l_4) \), then we have \( |[\{x^l_{l_1}, u_3\}, V(P^r_4)]_G| \geq 6 \). Therefore, the subgraph \( G[\{x^l_{l_1}, y^i_1, y^i_2, u_3\} \cup V(P^r_4)] \) contains a rainbow subgraph isomorphic to \( 2P_4 \).\\
\noindent
{\bf 1.3} If \( u_1u_2 \notin E(P^l_4) \) and \( x^l_{l_1}u_3 \in E(P^l_4) \), then it follows that \( |[\{u_1, u_2\}, V(P^r_4)]_G| \geq 6 \). Thus, for \( s_1 \neq s_2 \in [2] \), there exist \( z_1 \in [2] \) and \( z_2 \in \{3, 4\} \) such that \( u_{s_1} x^r_{z_1}, u_{s_2} x^r_{z_2} \in E(G) \).\\
\noindent
{\bf 1.4} If both \( u_1u_2 \) and \( x^l_{l_1}u_3 \notin E(P^l_4) \), then for any two distinct indices \( z_1, z_2 \in [4] \), we have \( |[\{x^l_{z_1}, x^l_{z_2}\}, V(P^r_4)]_G| \geq 6 \). Consequently, there are four distinct indices \( s_i \), where \( i \in [4] \), in the set \( [4] \) such that \( \{x^r_{s_1}x^l_{l_1}, x^r_{s_2}u_3, x^r_{s_3}x^r_{s_4}\} \subseteq E(G) \). Additionally, for an element \( a \in [2] \), we have \( [\{u_a\}, \{x^r_{s_3}, x^r_{s_4}\}]_G \neq \emptyset \).

From propositions (1.1) to (1.4), we conclude that Case 1 cannot occur; otherwise, we could construct a rainbow subgraph of \( K_n \) that is isomorphic to \( kP_4 \cup tP_2 \), which contradicts our assumption.\\
\noindent
{\bf Case 2.} Consider two elements \( P^l_4 \) and \( P^r_4 \) in the set \( \mathcal{F}_h \). In this case, it is true that \( u_1, u_2 \in V(P^l_4) \) and \( u_3 \in V(P^r_4) \). By applying Equation \eqref{seq1}, we can assert the following statements:\\
{\bf 2.1} If \( u_1 = x^l_1 \), \( u_2 = x^l_2 \), and \( u_3 \in \{ x^r_1, x^r_2 \} \), then we have \( [\{ x^l_3, x^l_4 \}, \{ x^r_3, x^r_4 \}]_G \neq \emptyset \). \\
{\bf 2.2} If \( u_1 = x^l_2 \), \( u_2 = x^l_3 \), and \( u_3 \in \{ x^r_1, x^r_2 \} \), then either \( |[\{ x^l_1, x^l_4 \}, \{ x^r_3, x^r_4 \}]_G| \geq 3 \) or \( \left( |[\{ x^l_1, x^l_4 \}, \{ x^r_3, x^r_4 \}]_G| = 2 \text{ and } |[\{ x^l_1, x^l_4 \}, \{ x^r_1, x^r_2 \}]_G| = 4 \right) \) holds. \\
{\bf2.3} If \( u_1 = x^l_1 \), \( u_2 = x^l_3 \), and \( u_3 \in \{ x^r_1, x^r_2 \} \), then either \( |[\{ x^l_2, x^l_4 \}, \{ x^r_3, x^r_4 \}]_G| \geq 3 \) or \( (x^l_1 x^r_3 \in E(G) \text{ and } |[\{ x^l_2, x^l_4 \}, \{ x^r_1, x^r_2 \}]_G| = 4) \) holds. \\
{\bf 2.4} If \( u_1 = x^l_1 \), \( u_2 = x^l_4 \), and \( u_3 \in \{ x^r_1, x^r_2 \} \), then either \( |[\{ x^l_2, x^l_3 \}, \{ x^r_3, x^r_4 \}]_G| \geq 3 \) or \( x^l_1 x^r_3 \in E(G) \) holds.

By analyzing statements (2.1) to (2.4) along with the symmetry of the path \( P_4 \), we can conclude that Case 2 cannot occur. If it did, we would be able to construct a rainbow subgraph of \( K_n \) that is isomorphic to \( kP_4 \cup tP_2 \), which contradicts our initial assumption.\\
{\bf Case 3.} Four two elements \( P^l_4 \) and \( P^r_4 \) in the set \( \mathcal{F}_h \), it is true that \( u_1, u_3 \in V(P^l_4) \) and \( u_2 \in V(P^r_4) \). By applying Equation \eqref{seq1}, we can confirm the following conditions:\\
{\bf 3.1} If \( u_1 = x^l_1 \), \( u_3 = x^l_2 \), and \( u_2= x^r_1\), then either $|[\{x^l_1,x^l_4\},V(P^r_4)]_G|\geq7$ or ($|[\{x^l_1,x^l_4\},V(P^r_4)]_G|=6$ and $|[\{x^l_2,x^l_3\},V(P^r_4)]_G|=8$) holds.\\
{\bf 3.2} If \( u_1 = x^l_1 \), \( u_3 = x^l_3 \), and \( u_2= x^r_1\), then either $|[\{x^l_1,x^l_2\},V(P^r_4)]_G|\geq7$, $|[\{x^l_1,x^l_4\},V(P^r_4)]_G|\geq7$, or ($|[\{x^l_1\},V(P^r_4)]_G|=2$ and $|[\{x^l_2,x^l_4\},V(P^r_4)]_G|=8$) holds.\\
{\bf 3.3} If \( u_1 = x^l_1 \), \( u_3 = x^l_4 \), and \( u_2= x^r_1\), then either $|[\{x^l_1,x^l_2\},V(P^r_4)]_G|\geq7$ or $|[\{x^l_4\},V(P^r_4)]_G|=4$ holds.\\
{\bf 3.4} If \( u_1 = x^l_2 \), \( u_3 = x^l_3 \), and \( u_2= x^r_1\), then either $|[\{x^l_1,x^l_2\},V(P^r_4)]_G|\geq7$ or ($|[\{x^l_1\},V(P^r_4)]_G|\geq2$ and $|[\{x^l_3,x^l_4\},V(P^r_4)]_G|=8$) holds.\\
{\bf 3.5} If \( u_1 = x^l_2 \), \( u_3 = x^l_4 \), and \( u_2= x^r_1\), then either $|[\{x^l_1,x^l_2\},V(P^r_4)]_G|\geq7$ or ($|[\{x^l_1\},V(P^r_4)]_G|\geq2$ and $|[\{x^l_3,x^l_4\},V(P^r_4)]_G|=8$) holds.\\
{\bf 3.6} If \( u_1 = x^l_3 \), \( u_3 = x^l_4 \), and \( u_2= x^r_1\), then either $x^l_4x^r_4\in E(G)$, ($x^l_4x^r_2\in E(G)$ and $|[\{x^l_1,x^l_2,x^l_3\},V(P^r_4)]_G|\geq11$), or ($x^l_4x^r_3\in E(G)$ and $|[\{x^l_1,x^l_2,x^l_3\},V(P^r_4)]_G|=12$) holds.\\
{\bf 3.7} If \( u_1 = x^l_1 \), \( u_3 = x^l_2 \), and \( u_2= x^r_2\), then either $x^l_1x^r_1,x^r_1x^l_4\in E(G)$, ($x^r_1x^l_2\in E(G)$ and $|[\{x^r_2,x^r_3,x^r_4\},V(P^l_4)]_G|\geq11$ ), or ($x^r_1x^l_3\in E(G)$ and $|[\{x^r_3,x^r_4\},V(P^l_4)]_G|=8$ ) holds.\\
{\bf 3.8} If \( u_1 = x^l_1 \), \( u_3 = x^l_3 \), and \( u_2= x^r_2\), then either $x^r_1x^l_2\in E(G)$, ($x^r_1x^l_3\in E(G)$ and $|[\{x^r_2,x^r_3,x^r_4\},V(P^l_4)]_G|\geq11$), or ($x^r_1x^l_4\in E(G)$ and $|[\{x^r_3,x^r_4\},V(P^l_4)]_G|=8$) holds.\\
{\bf 3.9} If \( u_1 = x^l_1 \), \( u_3 = x^l_4 \), and \( u_2= x^r_2\), then either $x^r_1x^l_2\in E(G)$, $x^r_1x^l_4\in E(G)$, or ($x^r_1x^l_3\in E(G)$ and $|[\{x^r_3,x^r_4\},V(P^l_4)]_G|=8$) holds.\\
{\bf 3.10} If \( u_1 = x^l_2 \), \( u_3 = x^l_3 \), and \( u_2= x^r_2\), then either $x^r_1x^l_1\in E(G)$, ($x^r_1x^l_3,x^r_2x^l_2\in E(G)$ and $|[\{x^l_1,x^l_4\},\{x^r_3,x^r_4\}]_G|\geq3$), 
($\{x^r_1x^l_3,x^r_1x^l_4\}\cap E(G)\neq\emptyset$ and $|[\{x^r_3,x^r_4\},V(P^l_4)]_G|=8$) holds. \\
{\bf 3.11} If \( u_1 = x^l_2 \), \( u_3 = x^l_4 \), and \( u_2= x^r_2\), then either $x^r_1x^l_1\in E(G)$, ($x^r_1x^l_4,x^r_2x^l_2\in E(G)$ and $|[\{x^l_1,x^l_3\},\{x^r_3,x^r_4\}]_G|\geq3$), or ($\{x^r_1x^l_3,x^r_1x^l_4\}\cap E(G)\neq\emptyset$ and $|[\{x^r_3,x^r_4\},V(P^l_4)]_G|=8$) holds. \\
{\bf 3.12} If \( u_1 = x^l_3 \), \( u_3 = x^l_4 \), and \( u_2= x^r_2\), then either $x^r_1x^l_4\in E(G)$, ($x^r_1x^l_1\in E(G)$ and $|[\{x^r_3,x^r_4\},\{x^l_2,x^l_3,x^l_4\}]_G|\geq5$), or ($x^r_1x^l_2\in E(G)$ and $|[\{x^r_3,x^r_4\},\{x^l_1,x^l_3,x^l_4\}]_G|=6$) holds.

By utilizing conditions (3.1) to (3.12) along with the symmetry of the path \( P_4 \), we can conclude that Case 3 cannot hold true. If it did, we would be able to construct a rainbow subgraph of \( K_n \) that is isomorphic to \( kP_4 \cup tP_2 \), which contradicts our initial assumption.\\
{\bf Case 4.} In this scenario, we have three elements \( P^a_4 \), \( P^b_4 \), and \( P^s_4 \) in the set \( \mathcal{F}_h \). It is given that \( u_1 \in V(P^a_4) \), \( u_2 \in V(P^b_4) \), and \( u_3 \in V(P^s_4) \). We can assume that \( x^a_{a_1} \) is a vertex of \( V(P^a_4) \) such that \( x^a_{a_1} \neq u_1 \), and that the subgraph \( G[(V(P^a_4)-\{x^a_{a_1}\})\cup\{v_1\}] \) contains a rainbow subgraph isomorphic to \( P_4 \). 
By applying Equation \eqref{seq1}, there exists a vertex \( x^b_{b_1} \) in \( V(P^b_4) \) such that \( x^b_{b_1} \neq u_2 \) and \( x^b_{b_1}x^a_{a_1} \in E(G) \). 
Therefore, if we replace \( v_2 \) and \( u_2 \) with \( x^a_{a_1} \) and \( x^b_{b_1} \), respectively, and also substitute \( P^l_4 \) and \( P^r_4 \) with \( P^b_4 \) and \( P^s_4 \), we can apply Case 2 to achieve our objective.\qed

\begin{claim}\label{scl2}
Let \(h \in [\lambda]\), $i\neq j\in[t+2]$, and $i_1,j_1\in[2]$. If $|\mathcal{F}_h|\geq2$, then there are not two distinct vertices $u_1$ and $u_2$ in $\mathcal{V}_h$ such that $\{y^i_{i_1}u_1,y^j_{j_1}u_2\}\subseteq E(G)$. 
\end{claim}
	
Assume, for the sake of contradiction, that \( |\mathcal{F}_h| \geq 2 \) and that there are two distinct vertices \( u_1 \) and \( u_2 \) in \( \mathcal{V}_h \) such that \( \{ y^i_{i_1} u_1, y^j_{j_1} u_2 \} \subseteq E(G) \). Based on the positions of the vertices \( u_1 \) and \( u_2 \), we will consider two cases as follows:\\
{\bf Case 1.} For an element \( P^l_4 \) in the set \( \mathcal{F}_h \), we have \( u_1, u_2 \in V(P^l_4) \). In this case, if there are two distinct indices \( l_1, l_2 \in [4] \) such that \( x^l_{l_1} u_1 \) and \( x^l_{l_2} u_2 \) are both edges in \( E(P^l_4) \), and if \( \{ x^l_{l_1}, x^l_{l_2} \} \cap \{ u_1, u_2 \} =\emptyset \), then for \( i_2 \in ([2] - \{i_1\}) \) and \( j_2 \in ([2] - \{j_1\}) \), the paths \( y^i_{i_2} y^i_{i_1} u_1 x^l_{l_1} \) and \( y^j_{j_2} y^j_{j_1} u_2 x^l_{l_2} \) create a rainbow subgraph of \( G \) that is isomorphic to \( 2P_4 \), which leads to a contradiction. 
	
Otherwise, there exist \( a \in \{ 1, 3 \} \) and \( s \neq z \in \{ i, j \} \) such that \( y^s_{s_1} x^l_{a}, y^z_{z_1} x^l_{a+1} \in E(G) \). In this scenario, assume \( P^r_4 \) is a path different from \( P^l_4 \) in the set \( \mathcal{F}_h \). Based on Equation (1), there are two distinct indices \( s \) and \( z \) in the set \( [4] \) such that \( x^r_s u_1 \) and \( x^r_z u_2 \) are edges in \( E(G) \), \( P^r_4 - \{ x^r_s, x^r_z \} \cong P_2 \), and \( [V(P^l_4)-\{ x^l_{a}, x^l_{a+1} \}, V(P^r_4)-\{ x^r_s, x^r_z \}]_G \neq \emptyset \). These conditions ensure that the subgraph \( G[V(P^i_2) \cup V(P^j_2) \cup V(P^l_4) \cup V(P^r_4)] \) contains a rainbow subgraph isomorphic to \( 3P_4 \), leading to a contradiction.\\
{\bf Case 2.} For two distinct elements \( P^l_4 \) and \( P^r_4 \) in the set \( \mathcal{F}_h \), we have \( u_1 \in V(P^l_4) \) and \( u_2 \in V(P^r_4) \). In this scenario, let \( s \) and \( z \) be two distinct indices in the set \([4]\) such that \( x^l_s u_1, x^r_z u_2 \in E(G) \). Additionally, we assume that \( P^l_4 - \{x^l_s, u_1\} \cong P_2 \) and \( P^r_4 - \{x^r_z, u_2\} \cong P_2 \). 
By applying Equation \eqref{seq1}, we can observe that the edge set \([V(P^l_4) - \{x^l_s, u_1\}, V(P^r_4) - \{x^r_z, u_2\}]_G\)  is not empty. This implies that the subgraph \( G[V(P^i_2) \cup V(P^j_2) \cup V(P^l_4) \cup V(P^r_4)] \) contains a rainbow subgraph that is isomorphic to \( 3P_4 \), which leads to a contradiction. \qed

\begin{claim}\label{scl3}
Let \(h \in [\lambda]\). If $|\mathcal{F}_h|\neq2$, then  there are not four distinct vertices $v_i$, $i\in[4]$, in $V'$ and four distinct vertices $u_i$, $i\in [4]$, in $\mathcal{V}_h$ such that $\{v_iu_i:i\in[4]\}\subseteq E(G)$. 
\end{claim}

To prove this claim, we will assume, for the sake of contradiction, that there are four distinct vertices \( v_i \) for \( i \in [4] \) in \( V' \) and four distinct vertices \( u_i \) for \( i \in [4] \) in \( \mathcal{V}_h \) such that the set \( \{v_i u_i : i \in [4]\} \) is a subset of \( E(G) \). Based on the positions of the vertices \( u_i \), we will consider the following cases:\\
{\bf Case 1.} For an element \( P^a_4 \) in \( \mathcal{F}_h \), we have \( \{u_i : i \in [4]\} = V(P^a_4) \). In this situation, it is clear that the subgraph \( G[V(P^a_4) \cup \{v_i : i \in [4]\}] \) contains a rainbow subgraph isomorphic to \( 2P_4 \). However, this contradicts our assumption that \( K_n \) does not have a rainbow subgraph isomorphic to \( kP_4 \cup tP_2 \). Therefore, this case cannot occur.\\
{\bf Case 2.} There are two distinct elements, \( P^a_4 \) and \( P^b_4 \), in \( \mathcal{F}_h \) such that for an element \( \{i_1, i_2\} \) in the set \( \{\{1, 2\}, \{3, 4\}\} \), it holds that 
\(V(P^a_4) \cap \{u_i : i \in [4]\} = \{x^a_{i_1}, x^a_{i_2}\} \)
and \( |V(P^b_4) \cap \{u_i : i \in [4]\}| \leq 1 \). Under these circumstances, by employing Equation \eqref{seq1}, there exists an index \( j \in [4] \) such that
 \[ x^b_j \not\in \{u_i : i \in [4]\}~{\rm and}~ [\{x^b_j\}, V(P^a_4) - \{x^a_{i_1}, x^a_{i_2}\}]_G \neq \emptyset. \] 
Therefore, using a method similar to that in the proof of Claim \ref{scl1}, we can conclude that Case 2 cannot occur. \\
{\bf Case 3.} There are two distinct elements, \( P^a_4 \) and \( P^b_4 \), in \( \mathcal{F}_h \) such that 
\[ V(P^a_4) \cap \{u_i : i \in [4]\} = \{x^a_{1}, x^a_{4}\} \quad \text{and} \quad |V(P^b_4) \cap \{u_i : i \in [4]\}| \leq 1. \]
In this case, by applying Equation \eqref{seq1}, we find an index \( j \in [4] \) such that 
$x^b_j \not\in \{u_i : i \in [4]\}$ and $x^b_jx^a_4\in E(G)$. Thus, by utilizing a method similar to the one in the proof of Claim \ref{scl1}, we can conclude that Case 3 is not possible.\\
{\bf Case 4.} There are three distinct elements, \( P^a_4 \), \( P^b_4 \), and \( P^s_4 \), in the set \( \mathcal{F}_h \) such that for two elements \(\{i_1, i_2\}\) and \(\{j_1, j_2\}\), which are not necessarily distinct, chosen from the set \(\{\{1, 3\}, \{2, 4\}\}\), we have $V(P^a_4)\cap \{u_i:i\in [4]\}=\{x^a_{i_1},x^a_{i_2}\}$ and $V(P^b_4)\cap \{u_i:i\in [4]\}=\{x^b_{j_1},x^b_{j_2}\}$. For simplicity, we can assume that \(\{i_1, i_2\} = \{j_1, j_2\} = \{1, 3\}\), \(u_1 = x^a_1\), \(u_2 = x^a_3\), \(u_3 = x^b_1\), and \(u_4 = x^b_3\). Given these assumptions, we can utilize Equation \eqref{seq1} to observe that at least one of the following conditions must hold true.\\
{\bf 4.1} For each \( z \in \{a, b\} \) and \( i \in \{2, 4\} \), it holds that \( |[V(P^s_4), \{x^z_i\}]_G| \geq 2 \). Additionally, there exists a two-element subset \( W \) of the set \( \{x^a_2, x^a_4, x^b_2, x^b_4\} \) such that \( |[V(P^s_4), W]_G| = 8 \).\\
{\bf 4.2} There exists a \( z \in \{a, b\} \) such that \( |[V(P^s_4), \{x^z_1, x^z_3\}]_G| = 8 \). Furthermore, either \( \{x^a_2 x^b_2, x^a_4 x^b_4\} \) or \( \{x^a_2 x^b_4, x^a_4 x^b_2\} \) is a subset of \( E(G) \).\\
{\bf 4.3} There are \( z_1 \in \{a, b\} \), \( z_2 \in (\{a, b\} - \{z_1\}) \), \( i \in \{2, 4\} \), and \( j \in (\{2, 4\} - \{i\}) \) such that \( |[\{x^{z_1}_i\}, \{x^{z_2}_2, x^{z_2}_4\}]_G| = 2 \) and \( |[\{x^{z_1}_j\}, \{x^{z_2}_2, x^{z_2}_4\}]_G| = 0 \). Moreover, one of the following sets is a subset of \( E(G) \): \( \{x^{z_1}_j x^s_2, x^{z_2}_2 x^s_1\} \), \( \{x^{z_1}_j x^s_2, x^{z_2}_4 x^s_1\} \), \( \{x^{z_1}_j x^s_3, x^{z_2}_2 x^s_4\} \), \( \{x^{z_1}_j x^s_3, x^{z_2}_4 x^s_4\} \), \( \{x^{z_1}_j x^s_1, x^{z_2}_2 x^s_4\} \), \( \{x^{z_1}_j x^s_1, x^{z_2}_4 x^s_4\} \), \( \{x^{z_1}_j x^s_4, x^{z_2}_2 x^s_1\} \), \( \{x^{z_1}_j x^s_4, x^{z_2}_4 x^s_1\} \), \(\{x^{z_1}_j x^s_1,$ $x^{z_1}_j x^s_3,x^{z_2}_2 x^s_2\}\), \(\{x^{z_1}_j x^s_1,$ $x^{z_1}_j x^s_3,x^{z_2}_4 x^s_2\}\), \(\{x^{z_1}_j x^s_2,x^{z_1}_j x^s_4,x^{z_2}_2 x^s_3\}\), or \(\{x^{z_1}_j x^s_2,x^{z_1}_j x^s_4,x^{z_2}_4 x^s_3\}\). Additionally, \( |[\{x^a_1, x^a_3\}, \{x^b_1, x^b_3\}]_G| = 4 \).

Utilizing the conditions from (4.1) to (4.3), we can conclude that Case 4 cannot hold. If it were true, we would be able to construct a rainbow subgraph of \( K_n \) that is isomorphic to \( kP_4 \cup tP_2 \), which contradicts our assumption.\\
{\bf Case 5.} There are three distinct elements, \( P^a_4 \), \( P^b_4 \), and \( P^s_4 \), in \( \mathcal{F}_h \) such that \( V(P^a_4) \cap \{u_i : i \in [4]\} = \{x^a_2, x^a_3\} \) and \( V(P^b_4) \cap \{u_i : i \in [4]\} = \{x^b_2, x^b_3\} \). Under these circumstances, using Equation \eqref{seq1}, we can observe that at least one of the following conditions holds:\\
{\bf 5.1} There exists a two-element subset \( W \) of the set \( \{x^a_1, x^a_4, x^b_1, x^b_4\} \) such that the cardinality of the set \([V(P^s_4), W]_G\) equals to \( 8 \). Furthermore, for each \( z \in \{a, b\} \) and \( i \in \{1, 4\} \), it holds that \( |[V(P^s_4), \{x^z_i\}]_G| \geq 2 \).\\
{\bf 5.2} There exists \( z \in \{a, b\} \) such that \( |[V(P^s_4), \{x^z_2, x^z_4\}]_G| \geq 8 \). Additionally, either \( \{x^a_1 x^b_1, x^a_4 x^b_4\} \) or \( \{x^a_1 x^b_4, x^a_4 x^b_1\} \) is a subset of \( E(G) \).\\
{\bf 5.3} There exist \( z_1 \in \{a, b\} \), \( z_2 \in (\{a, b\} - \{z_1\}) \), \( i \in \{1, 4\} \), and \( j \in (\{1, 4\} - \{i\}) \) such that \( |[\{x^{z_1}_i\}, \{x^{z_2}_1, x^{z_2}_4\}]_G| = 2 \) and \( |[\{x^{z_1}_j\}, \{x^{z_2}_1, x^{z_2}_4\}]_G| = 0 \). Moreover, at least one of the following sets is a subset of \( E(G) \): \( \{x^{z_1}_j x^s_2, x^{z_2}_1 x^s_1\} \), \( \{x^{z_1}_j x^s_2, x^{z_2}_4 x^s_1\} \), \( \{x^{z_1}_j x^s_3, x^{z_2}_1 x^s_4\} \), \( \{x^{z_1}_j x^s_3, x^{z_2}_4 x^s_4\} \), \( \{x^{z_1}_j x^s_1, x^{z_2}_1 x^s_4\} \), \( \{x^{z_1}_j x^s_1, x^{z_2}_4 x^s_4\} \), \( \{x^{z_1}_j x^s_4, x^{z_2}_1 x^s_1\} \), \( \{x^{z_1}_j x^s_4, x^{z_2}_4 x^s_1\} \), \( \{x^{z_1}_j x^s_1,\) \(x^{z_1}_j x^s_3, x^{z_2}_1 x^s_2\} \), \( \{x^{z_1}_j x^s_1, x^{z_1}_j x^s_3, x^{z_2}_4 x^s_2\} \), \( \{x^{z_1}_j x^s_2, x^{z_1}_j x^s_4, x^{z_2}_1 x^s_3\} \), and \( \{x^{z_1}_j x^s_2, x^{z_1}_j x^s_4, x^{z_2}_4 x^s_3\} \).

Utilizing the conditions from (5.1) to (5.3), we can conclude that Case 5 cannot hold. If it were true, we would be able to construct a rainbow subgraph of \( K_n \) that is isomorphic to \( kP_4 \cup tP_2 \), which contradicts our assumption.\\
{\bf Case 6.} There are three distinct elements, \( P^a_4 \), \( P^b_4 \), and \( P^s_4 \), in \( \mathcal{F}_h \) such that \( |V(P^a_4) \cap \{u_i : i \in [4]\}| = 3 \) and \( |V(P^b_4) \cap \{u_i : i \in [4]\}| = 1 \). Let \( i \) be an index in \( [3] \) such that \( x^a_i, x^a_{i+1} \in \{u_i : i \in [4]\} \). Define \( [4] - \{i, i+1\} = \{l_1, l_2\} \). Under these conditions, using Equation \eqref{seq1}, we observe that there are two distinct indices \( j_1 \) and \( j_2 \) in the set \( [4] \) such that \( x^a_{l_1} x^s_{j_1} \) and \( x^a_{l_2} x^s_{j_2} \) are in \( E(G) \). Therefore, we can apply a proof similar to that in Claim \ref{scl1} to demonstrate that Case 6 cannot occur.\\
{\bf Case 7.}  There are three distinct elements, \( P^a_4 \), \( P^b_4 \), and \( P^s_4 \), within \( \mathcal{F}_h \) such that \( V(P^a_4) \cap \{u_i : i \in [4]\} = \{x^a_1,x^a_3\} \) or \( \{x^a_2,x^a_4\} \), \( V(P^b_4) \cap \{u_i : i \in [4]\} = \{x^b_1\} \) or \( \{x^b_4\} \), and \( V(P^s_4) \cap \{u_i : i \in [4]\} = \{x^s_1\} \) or \( \{x^s_4\} \). Without loss of generality, we can assume that \( V(P^a_4) \cap \{u_i : i \in [4]\} = \{x^a_1,x^a_3\} \), \( V(P^b_4) \cap \{u_i : i \in [4]\} = \{x^b_1\} \), and \( V(P^s_4) \cap \{u_i : i \in [4]\} = \{x^s_1\} \). Given these assumptions, we can refer to Equation \eqref{seq1} to observe that at least one of the following conditions must hold true:\\
{\bf 7.1} There exist \( z_1 \neq z_2 \in \{b,s\} \) and \( i_1 \neq i_2 \in \{2,4\} \) such that either $\{x^{a}_{i_1}x^{z_1}_4,x^{a}_{i_2}x^{z_2}_4\}$, $\{x^{a}_{i_1}x^{z_1}_4,x^{a}_{i_2}x^{z_2}_2\}$, $\{x^{a}_{i_1}x^{z_1}_2,x^{a}_{i_2}x^{z_2}_4\}$, or $\{x^{a}_{i_1}x^{z_1}_2,x^{a}_{i_2}x^{z_2}_2\}$ is a subset of $E(G)$.
Furthermore, the sizes of each of the sets \( [V(P^{b}_4),\{x^{a}_{1},x^{a}_{2}\}]_G \) and \( [V(P^{s}_4),\{x^{a}_{3},x^{a}_{4}\}]_G \) are both greater than or equal to $6$.\\
{\bf 7.2} There exists \( i \in \{2,4\} \) such that \( [V(P^{b}_4), \{x^{a}_{i}\}]_G = \{x^{a}_{i}x^b_1, x^{a}_{i}x^b_3\} \) and \( [V(P^{s}_4), \{x^{a}_{i}\}]_G$ $=$ $ \{x^{a}_{i}x^s_1, x^{a}_{i}x^s_3\} \). Additionally, either \( x^b_2x^s_2 \), \( x^b_2x^s_4 \), or \( x^b_4x^s_2 \) is an edge in \( E(G) \). It is also true that \( [\{x^b_1\}, V(P^s_4) - \{x^s_{2}\}]_G \neq \emptyset \).

By examining the conditions outlined in (7.1) and (7.2), we can conclude that Case 7 cannot hold. If it were true, we would be able to construct a rainbow subgraph of \( K_n \) that is isomorphic to \( kP_4 \cup tP_2 \), which contradicts our initial assumption.\\
{\bf Case 8.}  There are three distinct elements, \( P^a_4 \), \( P^b_4 \), and \( P^s_4 \), within \( \mathcal{F}_h \) such that \( V(P^a_4) \cap \{u_i : i \in [4]\} = \{x^a_1, x^a_3\} \) or \( \{x^a_2, x^a_4\} \), \( V(P^b_4) \cap \{u_i : i \in [4]\} = \{x^b_2\} \) or \( \{x^b_3\} \), and  \( V(P^s_4) \cap \{u_i : i \in [4]\} = \{x^s_2\} \) or \( \{x^s_3\} \).
Without loss of generality, we can assume \( V(P^a_4) \cap \{u_i : i \in [4]\} = \{x^a_1, x^a_3\} \), \( V(P^b_4) \cap \{u_i : i \in [4]\} = \{x^b_2\} \), and  \( V(P^s_4) \cap \{u_i : i \in [4]\} = \{x^s_2\} \). Under these assumptions, we refer to Equation \eqref{seq1}, which indicates that at least one of the following conditions must be satisfied:\\
{\bf 8.1} There exist \( z_1 \neq z_2 \in \{b,s\} \) and \( i_1 \neq i_2 \in \{2,4\} \) such that \( \{x^{a}_{i_1}x^{z_1}_1, x^{a}_{i_2}x^{z_2}_1\} \) forms a subset of \( E(G) \).\\
{\bf 8.2} There is a \( z \in \{b,s\} \) such that \( |[\{x^a_2,x^a_4\},\{x^z_3,x^z_4\}]_G| \geq 3 \). Additionally, \( x^b_1 x^s_1 \in E(G) \).\\
{\bf 8.3} There are \( z_1 \in \{b,s\} \), \( z_2 \in (\{b,s\} - \{z_1\}) \), and \( i_1 \neq i_2 \in \{3,4\} \) such that
\[
\{x^{a}_{2}x^{z_1}_{1}, x^{a}_{4}x^{z_1}_{1}, x^{a}_{4}x^{z_2}_{4}, x^{z_2}_1x^{z_1}_{i_1}, x^{z_2}_3x^{z_1}_{i_2}\}\subseteq E(G).
\]
\noindent
{\bf 8.4} There are \( z_1 \neq z_2 \in \{b,s\} \), \( i_1 \in \{1,3\} \), \( i_2 \in [4] \), and \( j_1 \in \{2,4\} \) such that \( x^a_{i_1} x^{z_1}_{i_2} \notin E(G) \) and \( x^a_{j_1} x^{z_1}_1 \in E(G) \). Additionally, \( |[\{x^a_2,x^a_4\}, \{x^{z_1}_3, x^{z_1}_4\}]_G| = 4 \), and either \( x^{z_1}_3 x^{z_2}_1 \in E(G) \) or \( x^{z_1}_4 x^{z_2}_1,x^{z_1}_2 x^{z_2}_1 \in E(G) \) holds.\\
{\bf 8.5} There are \( z_1 \neq z_2 \in \{b,s\} \), \( i_1, i_2 \in \{2,4\} \), and \( j_1 \neq j_2 \in \{3,4\} \), such that \( \{x^a_{i_1}x^{z_1}_1, x^a_{i_2}x^{z_1}_{j_1}, x^{z_1}_{j_2}x^{z_2}_1\} \subseteq E(G) \). Additionally, \( |[\{x^a_1,x^a_3\},\{x^{z_1}_2\}]_G| = 2 \).\\
{\bf 8.6} There are \( z_1 \neq z_2 \in \{b,s\} \), \( i_1, i_2 \in \{2,4\} \), and \( j_1 \neq j_2, l_1 \neq l_2 \in \{3,4\} \) such that
$\{x^a_{i_1}x^{z_1}_1, x^a_{i_1}x^{z_1}_{j_1}, x^{z_1}_{j_1}x^{z_2}_1, x^a_{i_2}x^{z_2}_{l_1}, x^{z_1}_{j_2}x^{z_2}_{l_2}\} \subseteq E(G)$.
Additionally, \(\{x^a_1x^{z_1}_2,x^a_3x^{z_1}_2\}\subseteq E(G)\) and  \(|[\{x^{z_1}_2, x^{z_1}_3\}, \{x^{z_2}_2, x^{z_2}_3,x^{z_2}_4\}]_G|=6 \).\\
{\bf 8.7} All three of the conditions \( |[\{x^a_2,x^a_4\},\{x^{b}_3,x^b_4,x^s_3,x^s_4\}]_G| = 8 \), \( |[\{x^{b}_1\},\{x^{s}_3,x^{s}_4\}]_G| \geq 1 \), and \( |[\{x^{s}_1\},\{x^{b}_3,x^{b}_4\}]_G| \geq 1 \) holds.

By analyzing the conditions outlined in equations (8.1) and (8.7), we can conclude that Case 8 cannot be valid. If it were true, we would be able to construct a rainbow subgraph of \( K_n \) that is isomorphic to \( kP_4 \cup tP_2 \). This scenario contradicts our initial assumptions.\\
{\bf Case 9.} There are three distinct elements, \( P^a_4 \), \( P^b_4 \), and \( P^s_4 \), within \( \mathcal{F}_h \) such that \( V(P^a_4) \cap \{u_i : i \in [4]\} = \{x^a_1, x^a_3\} \) or \( \{x^a_2, x^a_4\} \), \( V(P^b_4) \cap \{u_i : i \in [4]\} = \{x^b_1\} \) or \( \{x^b_4\} \), and \( V(P^s_4) \cap \{u_i : i \in [4]\} = \{x^s_2\} \) or \( \{x^s_3\} \).
Without loss of generality, we can assume that \( V(P^a_4) \cap \{u_i : i \in [4]\} = \{x^a_1, x^a_3\} \), \( V(P^b_4) \cap \{u_i : i \in [4]\} = \{x^b_1\} \), and  \( V(P^s_4) \cap \{u_i : i \in [4]\} = \{x^s_2\} \). Under these assumptions, we refer to Equation \eqref{seq1}, which indicates that at least one of the following conditions must be satisfied:\\
{\bf 9.1} There exist \( i_1 \neq i_2 \in \{2,4\} \) such that \( \{x^{a}_{i_1}x^{b}_4, x^{a}_{i_2}x^{s}_1\} \) forms a subset of \( E(G) \).\\
{\bf 9.2} There are \( i_1 \neq i_2 \in \{2,4\} \) such that \( \{x^b_2x^s_1,x^{a}_{i_1}x^{b}_3, x^{a}_{i_2}x^{b}_4\} \) forms a subset of \( E(G) \).\\
{\bf 9.3} There are \( i_1 \neq i_2 \in \{2,4\} \) such that \( \{x^b_4x^s_1,x^{a}_{i_1}x^{s}_3, x^{a}_{i_2}x^{s}_4\} \) forms a subset of \( E(G) \).\\
{\bf 9.4} Both conditions $x^b_4x^s_1\in E(G)$ and  $[\{x^a_2,x^a_4\},\{x^b_4\}]_G=\emptyset$ hold. Additionally, $|[V(P^a_4),V(P^b_4)-\{x^b_4\}]_G|=12$.\\
{\bf 9.5} There are \( i_1 \neq i_2 \in \{2,4\} \) such that \( \{x^{a}_{i_1}x^{s}_3, x^{a}_{i_2}x^{s}_4,x^b_2x^s_1 \} \) forms a subset of \( E(G) \). Additionally, either $x^b_1x^s_2\in E(G)$ or  $[\{x^b_3,x^b_4\},\{x^s_2\}]_G\neq\emptyset$ holds.\\
{\bf 9.6} Both conditions $\{x^a_2x^s_1,x^a_4x^s_1,x^b_2x^s_1\}\subseteq E(G)$ and $|[\{x^a_2,x^a_4\},\{x^b_3,x^b_4\}]_G|=2$ hold.  Additionally, $[\{x^b_3,x^b_4\},\{x^s_2\}]_G\neq\emptyset$.\\ 
{\bf 9.7} There are \( i_1 \neq i_2 \in \{2,4\} \) and $j_1\neq j_2\in\{3,4\}$ such that $\{x^a_{i_1}x^b_4,x^a_{i_2}x^s_{j_1},x^b_2x^s_{j_2},x^b_2x^s_1\}$ forms a subset of $E(G)$. Additionally, either $\{x^b_3x^s_1\}$ or $\{x^b_3x^s_2x^b_3x^s_{j_2}\}$ forms a subset of $E(G)$.\\
{\bf 9.8} There are \( i_1 \neq i_2 \in \{2,4\} \) such that $\{x^a_{i_1}x^s_1,x^a_{i_2}x^b_{3},x^b_4x^s_{3},x^b_2x^s_4\}$ forms a subset of $E(G)$.\\
{\bf 9.10} The set $\{x^a_2x^b_2,x^a_4x^s_4,x^b_3x^s_1,x^b_4x^s_3,x^b_1x^s_2\}$ forms a subset of $E(G)$.

By examining the conditions outlined in (9.1) and (9.10), we can conclude that Case 9 cannot hold. If it were true, we would be able to construct a rainbow subgraph of \( K_n \) that is isomorphic to \( kP_4 \cup tP_2 \), which contradicts our initial assumptions.\\
{\bf Case 10.} There are four distinct elements, \( P^a_4 \), \( P^b_4 \), \( P^d_4 \), and \( P^s_4 \), within \( \mathcal{F}_h \). For each \( z \in \{a, b, d, s\} \), we have \( |V(P^z_4) \cap \{u_i : i \in [4]\}| = 1 \). Without loss of generality, we can assume that \( V(P^a_4) \cap \{u_i : i \in [4]\} = \{u_1\} \).
By using Equation \eqref{seq1}, we can observe that there exists an index \( i \in [4] \) such that \( G[(V(P^{a}_4) - \{x^{a}_{i}\}) \cup \{v_1\}] \) contains a rainbow subgraph that is isomorphic to \( P_4 \). Additionally, there exists an index \( j \in [4] \) such that \( x^b_j \not \in (V(P^b_4) \cap \{u_i : i \in [4]\}) \), and it also holds that \( x^a_{i}x^b_j \in E(G) \). Consequently, we can apply a method similar to that used in Case 9 to demonstrate that Case 10 cannot occur. \qed

\begin{claim}\label{scl4}
Let \(h \in [\lambda]\). If $|\mathcal{F}_h|=2$, then  there are not six distinct vertices $v_i$, $i\in[6]$, in $V'$ and six distinct vertices $u_i$, $i\in [6]$, in $\mathcal{V}_h$ such that $\{v_iu_i:i\in[6]\}\subseteq E(G)$. 
\end{claim}

To prove this claim, assume that $h$ is an index in $[\lambda]$ such that $|\mathcal{F}_h|=2$. Suppose, for $a, b \in [k-1]$, we have $\mathcal{F}_h = \{P_4^a, P_4^b\}$. For the sake of contradiction, assume there exist six distinct vertices $v_1, \ldots, v_6$ in $V'$, and six distinct vertices $u_1, \ldots, u_6$ in $\mathcal{V}_h$, such that $\{v_i u_i : i \in [6]\} \subseteq E(G)$. We analyze two cases:\\
{\bf Case 1.} There exists $z \in \{a, b\}$ such that $V(P_4^z)\subseteq\{u_i : i \in [4]\}$. In this situation, observe that the subgraph $G[V(P_4^z) \cup \{v_i : i \in [4]\}]$ contains a rainbow subgraph isomorphic to $2P_4$. However, this contradicts our assumption that $K_n$ does not contain a rainbow subgraph isomorphic to $kP_4 \cup tP_2$. Therefore, this case cannot occur.\\
{\bf Case 2.} For each $z \in \{a, b\}$, we have $|V(P_4^z) \cap \{u_i : i \in [6]\}| = 3$. In this case, by employing Equation~\eqref{seq1}, we observe that for each $i, j \in [3]$, $[\{x^a_i, x^a_{i+1}\}, \{x^b_j, x^b_{j+1}\}]_G \neq \emptyset$. This implies that the graph $G[V(P_4^a) \cup V(P_4^b) \cup \{v_i : i \in [6]\}]$ contains a rainbow subgraph isomorphic to $3P_4$. However, this again contradicts our initial assumption. \qed

Before examining the next claims, for \(h \in [\lambda]\), let’s define the set \(\mathcal{A}_h\) as follows:
\[\mathcal{A}_h = [\mathcal{V}_h, V(H_2) \cup V']_G.\]

\begin{claim}\label{scl5}
For \(h \in [\lambda]\), if $|\mathcal{F}_h|=1$, then
\[|\mathcal{A}_h| \leq\begin{cases}
   \max\{ 3|V'|, 2(|V'|+t+3) \} & {\rm when}~|V'|\geq4,\\
   \max\{ 4|V'|, 2(|V'|+t+3)\} & {\rm when}~|V'|\leq3.
  \end{cases}
\]
\end{claim}

Assume \( h \) is an index in the set \( [\lambda] \) such that \( |\mathcal{F}_h|=1 \). Additionally, for \( a\in [k-1] \), assume \( \mathcal{F}_h=\{P^a_4\} \). To demonstrate this claim, we will examine three cases based on the configurations of the edges in \( [V(P^a_4), V(H_2)]_G \) within \( P^a_4 \).\\
{\bf Case 1.} There exist indices \(i \in [t+2]\), \(i_1 \in [2]\), and \(l_1 \in \{1,4\}\) such that \(x^a_{l_1}y^i_{i_1} \in E(G)\). Without loss of generality, assume \(l_1 = 1\). Since \(K_n\) does not contain a rainbow subgraph isomorphic to \(kP_4 \cup tP_2\), we have the following statements:\\
{\bf 1.1} There are no \(j \in ([t+2]-\{i\})\), \(j_1 \in [2]\), and \(l_2 \in \{3,4\}\) such that \(x^a_{l_2}y^j_{j_1} \in E(G)\).\\
{\bf 1.2} There are no \(j \in ([t+2]-\{i\})\), \(j_1 \neq j_2 \in [2]\), and \(l_2 \in \{2,4\}\) such that both \(x^a_{1}y^j_{j_1}\) and \(x^a_{l_2}y^j_{j_2}\) are in \(E(G)\).\\
{\bf 1.3} There are no \(j \in ([t+2]-\{i\})\), \(j_1 \in [2]\), \(v \in V'\), and \(l_2 \in \{2,4\}\) such that either \(\{x^a_{1}v, x^a_{l_2}y^j_{j_1}\}\) or \(\{x^a_{1}y^j_{j_1}, x^a_{l_2}v\}\) forms a subset of \(E(G)\).\\
{\bf 1.4} There are no distinct vertices \(v_1\) and \(v_2\) in \(V'\) and \(l_2 \in \{2,4\}\) such that \(\{x^a_{1}v_1, x^a_{l_2}v_2\}\) forms a subset of \(E(G)\).\\
{\bf 1.5} There are no distinct vertices \(v_1\) and \(v_2\) in \(V'\) such that \(\{x^a_{3}v_1, x^a_{4}v_2\}\) forms a subset of \(E(G)\).

Based on the conditions established in statements (1.1) to (1.5), we can conclude the following: If for \( v \in (V' \cup (V(H_2) - \{y^i_1, y^i_2\})) \) it holds that \( vx^a_1 \in E(G) \), then we have \( |\mathcal{A}_h| \leq 2(|V'| + t + 2) \). Conversely, if this condition does not hold, then \( |\mathcal{A}_h| \leq 2(|V'| + t + 3) \).\\
{\bf Case 2.} There exist indices \(i \in [t+2]\), \(i_1 \in [2]\), and \(l_1 \in\{2,3\}\) such that \(x^a_{l_1}y^i_{i_1} \in E(G)\). Additionally, we have that \([\{x^a_1,x^a_4\}, V(H_2)]_G = \emptyset\). Since \(K_n\) does not contain a rainbow subgraph isomorphic to \(kP_4 \cup tP_2\), we can draw the following conclusions:\\
{\bf 2.1} There are no \(j \in ([t+2] - \{i\})\) and \(j_1 \in [2]\) such that \(x^a_{3}y^j_{j_1} \in E(G)\).\\
{\bf 2.2} There are no distinct vertices \(v_1 \neq v_2 \in V'\) and \(i_2 \in ([2] - \{i_1\})\) such that \([\{y^i_{i_2}, v_1, v_2\}, \{x^a_1, x^a_3, x^a_4\}]_G\) contains three independent edges.\\
{\bf 2.3} There are no distinct vertices \(v_1 \neq v_2 \in V'\) such that \(\{x^a_3v_1, x^a_4v_2\}\) is a subset of \(E(G)\).\\
{\bf 2.4} There are no distinct vertices \(v_1 \neq v_2 \in V'\) such that either \(\{x^a_1v_1,x^a_3v_1, x^a_2v_2\}\) or \(\{x^a_1v_1,x^a_4v_1, x^a_2v_2\}\) is a subset of \(E(G)\).\\
{\bf 2.5} There are no $v\in V'$, \(j \in ([t+2]-\{i\})\), and \(j_1 \in [2]\), such that either \(\{x^a_1v_1,x^a_3v_1, x^a_2y^j_{j_1}\}\) or \(\{x^a_1v_1,x^a_4v_1, x^a_2y^j_{j_1}\}\) is a subset of \(E(G)\).

From the conditions (2.1) to (2.3), we can observe that in this case, the cardinality of the set \( \mathcal{A}_h\) is at most \(2(|V'| + t + 2) \).\\
 {\bf Case 3.} There are no indices \(i \in [t+2]\), \(i_1 \in [2]\), and \(l_1 \in [4]\) such that \(x^a_{l_1}y^i_{i_1} \in E(G)\). Under this condition, let \(M\) represent the maximum matching in \(G\) where each edge has one endpoint in \(V'\) and the other in \(V(P^a_4)\). Given the maximality of \(M\) and our assumption that \(K_n\) does not contain a rainbow subgraph isomorphic to \(kP_4 \cup tP_2\), we can draw the following statements:\\
{\bf 3.1} The set \([V'-V(M), V(P^a_4)-V(M)]_G\) is empty, and the cardinality of the set \(M\) is at most \(3\).\\
{\bf 3.2} There are no vertices \(v\) in \((V'-V(M))\) and \(u\) in \((V(P^a_4)-V(M))\) such that for \(v_1\in (V'\cap V(M))\) and \(u_1 \in (V(P^a_4)\cap V(M))\), the edge \(v_1u_1\) belongs to \(M\) and both edges \(vu_1\) and \(uv_1\) are in \(E(G)\).

By applying statements (3.1) and (3.2), we conclude that \( |\mathcal{A}_h| \leq 3|V'| \) when \( |V'| \geq 4 \), and \( |\mathcal{A}_h| \leq 4|V'| \) when \( |V'| \leq 3 \), as required. \qed

\begin{claim}\label{scl6}
For \(h \in [\lambda]\), if $|\mathcal{F}_h|=2$, then
\[|\mathcal{A}_h| \leq
\begin{cases}
  \max\{5|V'|, 2(|V'| + t + 2)\} & {\rm when}~|V'|\geq8,\\
  \max\{40, 2(|V'| + t + 2)\} & {\rm when}~5\leq|V'|\leq7,\\
  \max\{8|V'|, 2(|V'| + t + 2)\} & {\rm when}~|V'|=4,\\
  \max\{3|\mathcal{V}_h|, 2(|V'| + t + 2)\} & {\rm when}~|V'|\leq3.   
\end{cases}
\]
\end{claim}

Let's suppose that for an index \( h \) in \( [\lambda] \), we have \( |\mathcal{F}_h| = 2 \). Additionally, assume that for \( a \neq b \in [k-1] \), it holds that \( \mathcal{F}_h = \{P^a_4, P^b_4\} \). To support this claim, we need to examine the following two cases:\\
{\bf Case 1.}  The set \( [\mathcal{V}_h, V(H_2)]_G \) is not empty. In this case, by employing Claim \ref{scl2}, we can observe that there are no \( i \neq j \in [t+2] \), \( i_1, j_1 \in [2] \), and \( u_1 \neq u_2 \in \mathcal{V}_h \) such that \( \{y^i_{i_1}u_1, y^j_{j_1}u_2\} \) forms a subset of \( E(G) \). Furthermore, applying Claim \ref{scl2} again confirms that there are no \( i \in [t+2] \), \( i_1 \in [2] \), \( v_1 \neq v_2 \in V' \), and \( u_1 \neq u_2 \neq u_3 \in \mathcal{V}_h \) such that \( \{y^i_{i_1}u_1, v_1u_2, v_2u_3\} \) forms a subset of \( E(G) \). Therefore, under these conditions, the cardinality of the set \( \mathcal{A}_h \) is at most \( \max\{3|\mathcal{V}_h|, 2|\mathcal{V}_h| + |V'|, 2(|V'| + t + 2)\} \).\\
{\bf Case 2.} The set \( [\mathcal{V}_h, V(H_2)]_G \) is empty. In this case, let $M$ denote a maximum matching in $G$ such that every edge in $M$ has one endpoint in $V'$ and the other in $\mathcal{V}_h$. By the maximality of $M$ and the assertions of Claim~\ref{scl4}, we observe the following:\\
{\bf 2.1} The size of $M$ is at most $5$. Furthermore, there are no vertices $u \in (V' - V(M))$ and $v \in (\mathcal{V}_h -V(M))$ such that $uv \in E(G)$.\\
{\bf 2.2} There do not exist vertices $v \in (V' - V(M))$ and $u \in (\mathcal{V}_h - V(M))$ for which, if $v_1 \in (V' \cap V(M))$ and $u_1 \in (\mathcal{V}_h  \cap V(M))$ with $v_1u_1 \in M$, both $vu_1$ and $uv_1$ belong to $E(G)$.

From statements (2.1) and (2.2), we deduce the following bounds:
If $|V'| \geq 8$, then $|\mathcal{A}_h| \leq 5|V'|$;
If $5 \leq |V'| \leq 7$, then $|\mathcal{A}_h| \leq 40$;
If $|V'| \leq 4$, then $|\mathcal{A}_h| \leq 8|V'|$.
This completes the argument for this case. \qed

\begin{claim}\label{scl7}
For \(h \in [\lambda]\), if $|\mathcal{F}_h|\geq3$, then
\[|\mathcal{A}_h| \leq
\begin{cases}
 \max\{3|V'|, 2(|V'| + t + 2)\} & {\rm when}~|V'|\geq|\mathcal{V}_h|,\\
 \max\{3|\mathcal{V}_h|, 2(|V'| + t + 2)\} & {\rm when}~|V'| \leq |\mathcal{V}_h|.
\end{cases}
\]
\end{claim}

Suppose there exists an index \( h \in [\lambda] \) such that \( |\mathcal{F}_h| \geq 3 \). To establish the claim, we consider the following two cases:\\
{\bf Case 1.} The set \( [\mathcal{V}_h, V(H_2)]_G \) is non-empty. In this scenario, applying an approach analogous to that used in Case 1 of the proof of Claim~\ref{scl6}, we find that the cardinality of \( \mathcal{A}_h \) is at most \( \max\{3|\mathcal{V}_h|,\ 2|\mathcal{V}_h| + |V'|,\ 2(|V'| + t + 2)\} \).\\
{\bf Case 2.} The set \( [\mathcal{V}_h, V(H_2)]_G \) is empty. Here, let $M$ be a maximum matching in $G$ in which each edge has one endpoint in $V'$ and the other in $\mathcal{V}_h$. By the maximality of $M$ and the conclusions of Claim~\ref{scl3}, we note the following:\\
{\bf 2.1} The size of $M$ is at most $3$. Moreover, there do not exist vertices $u \in (V' - V(M))$ and $v \in (\mathcal{V}_h - V(M))$ such that $uv \in E(G)$.\\
{\bf 2.2} There do not exist vertices $v \in (V' - V(M))$ and $u \in (\mathcal{V}_h - V(M))$ such that, for $v_1 \in (V' \cap V(M))$ and $u_1 \in (\mathcal{V}_h \cap V(M))$ with $v_1u_1 \in M$, both $vu_1$ and $uv_1$ belong to $E(G)$.

Consequently, from (2.1) and (2.2), we obtain the following bounds: If $|V'| \geq |\mathcal{V}_h|$, then $|\mathcal{A}_h| \leq 3|V'|$; If $3 \leq |V'| \leq |\mathcal{V}_h|$, then $|\mathcal{A}_h| \leq 3|\mathcal{V}_h|$; If $|V'| \leq 2$, then $|\mathcal{A}_h| \leq |V'||\mathcal{V}_h|$. This completes the argument for this case. \qed

Before proceeding to the final two claims, we introduce some additional notation. For each $i \in [\lambda]$, define
\[ \alpha_i = \begin{cases} \max\{3|V'|,\ 2(|V'| + t + 2)\} & \text{if } |\mathcal{F}_i| \geq 3,\\
 \max\{5|V'|,\ 2(|V'| + t + 2)\} & \text{if } |\mathcal{F}_i| = 2, \\
  \max\{3|V'|,\ 2(|V'| + t + 3)\} & \text{if } |\mathcal{F}i| = 1. \end{cases} \] 
Additionally, define
  \[\mathcal{S}=(|\mathcal{F}_1|,|\mathcal{F}_2|,\ldots,|\mathcal{F}_\lambda|)~{\rm and}~ \Omega(\mathcal{S})= \sum_{i=1}^{\lambda} \alpha_i+3\sum_{i=1}^{\lambda}|\mathcal{F}_i|(i-1). \]

\begin{claim}\label{scl8} 
If $|V'| \geq |V(H_1)|$, then $\Omega(\mathcal{S})$ is maximized when $\mathcal{S}=(\overbrace{1,1,\ldots,1}^{(k-1)~{\rm times}})$.
 \end{claim}
 To prove this claim, let $\mathcal{S}\neq(\overbrace{1,1,\ldots,1}^{(k-1)~{\rm times}})$ and $\alpha$ be the maximal index in $[\lambda]$ such that $|\mathcal{F}_\alpha| \neq 1$. We consider the proof in the following two cases:\\
{\bf Case 1.} $|\mathcal{F}_{\alpha}| = 2$.
Define $\mathcal{S}'=(|\mathcal{F}_1|,\ldots,|\mathcal{F}_{\alpha-1}|,\overbrace{1,\ldots,1}^{\upsilon~{\rm times}}),$ where $\upsilon = \lambda - \alpha + 2$. Using the definition of $\Omega$, we have \[ \Omega(\mathcal{S}') - \Omega(\mathcal{S}) = 2\max\{3|V'|, 2(|V'| + t + 3)\} - \max\{5|V'|, 2(|V'| + t + 2)\} - 3(\upsilon-1). \] 
Note that $\upsilon - 1 = \lambda - \alpha + 1 \leq k-2=\frac{|V(H_1)|-4}{4} \leq \frac{|V'|}{4}-1$. Additionally, one can observe that $\max\{6|V'|, 4(|V'| + t + 3)\} - \max\{5|V'|, 2(|V'| + t + 2)\}\geq |V'|$. Therefore, we can calculate that $\Omega(\mathcal{S}') - \Omega(\mathcal{S}) \geq\frac{1}{4}|V'|+3$. So, $\Omega(\mathcal{S}') > \Omega(\mathcal{S})$ in this case.\\
{\bf Case 2.} $|\mathcal{F}_{\alpha}| \geq 3$.
Now, set $\mathcal{S}'=(|\mathcal{F}_1|,\ldots,|\mathcal{F}_{\alpha-1}|,\overbrace{1,\ldots,1}^{\iota~{\rm times}})$, where $\iota = \lambda - \alpha + |\mathcal{F}_{\alpha}|$. Applying the definition of $\Omega$ yields 
\begin{align*}
\Omega(\mathcal{S}')-\Omega(\mathcal{S})&=|\mathcal{F}_{\alpha}|\max\{3|V'|,2(|V'| + t + 3)\}-\max\{3|V'|,2(|V'| + t + 2)\}\\
&~~~~~~~~~~~~~~~~~~~~~~~~~-3((\lambda-\alpha)(|\mathcal{F}_{\alpha}|-1)+\frac{1}{2}|\mathcal{F}_{\alpha}|(|\mathcal{F}_{\alpha}|-1))\\
&\geq(|\mathcal{F}_{\alpha}|-1)\left(\max\{3|V'|,2(|V'| + t + 3)\}-3(\lambda-\alpha)-\frac{3}{2}|\mathcal{F}_{\alpha}|\right).
\end{align*}
Since 
\[\lambda-\alpha+\frac{1}{2}|\mathcal{F}_{\alpha}|\leq k-1-|\mathcal{F}_{\alpha}|+1+\frac{1}{2}|\mathcal{F}_{\alpha}|= k-\frac{3}{2}=\frac{1}{4}|V(H_1)|-\frac{1}{2}\leq\frac{|V'|}{4}-\frac{1}{2},\]
we obtain
\[ \Omega(\mathcal{S}') - \Omega(\mathcal{S}) \geq (|\mathcal{F}_{\alpha}|-1)\left(\max\{3|V'|, 2(|V'| + t + 3)\} - \frac{3}{4}|V'|+\frac{3}{2}\right)>0. \] 
Thus, we also have  $\Omega(\mathcal{S}') > \Omega(\mathcal{S})$ in this case.

In both cases, we see that replacing $|\mathcal{F}_{\alpha}| > 1$ with ones increases $\Omega$. Therefore, $\Omega(\mathcal{S})$ is maximized when $\mathcal{S}=(\overbrace{1,1,\ldots,1}^{(k-1)~{\rm times}})$.\qed

By combining Claims~\ref{scl1} through~\ref{scl8}, along with some straightforward calculations, we can now derive the following upper bound for the size of the set $\mathcal{Y}(G)$:

\begin{claim}\label{scl9}
If $|V'| \geq |V(H_1)|$, then \[ |\mathcal{Y}(G)| \leq \frac{1}{2}(13k-14)(k-1) + (k-1) \max\{3|V'|,\ 2(|V'| + t + 3)\}. \] 
\end{claim}
 
We are now ready to complete the proof of our theorem. To this end, by combining Claims~\ref{scl0} and~\ref{scl9}, together with some straightforward calculations, we obtain the following upper bound on the size of $E(G)$:
\begin{equation}\label{seq2}
|E(G)|\leq\begin{cases}
  \frac{1}{2}(6kn+ 6t + 18 -11k^2 - 12kt - 19k - 2n )&{\rm if}~{n\geq4k+3t+3},\\
  \frac{1}{2}( 4kn+ 12 -3k^2 - 6kt - 21k )&{\rm if}~{n\leq4k+3t+3}.
\end{cases}
\end{equation}

On the other hand, by our assumption, we have \[|E(G)| = |c(K_n)| = AR(n, (2k + t)P_2) + 1.\] Since \(n \geq 8k + 2t - 4\) and \(k\geq t\geq2\), applying Equation~\eqref{eqth0} yields the following inequality:

\begin{equation}\label{seq3}
|E(G)|= (2k + t - 2)\left(n - \frac{2k+t-1}{2}\right) + 2.
\end{equation}

We define the following:
\begin{align*}
\beta_1 & = \frac{1}{2}(6kn+ 6t + 18 -11k^2 - 12kt - 19k - 2n ), \\
\beta_2 & = \frac{1}{2}( 4kn+ 12 -3k^2 - 6kt - 21k ), \\
\mu&= (2k + t - 2)\left(n - \frac{2k+t-1}{2}\right) + 2.
\end{align*} 

Let $f(x, y) = \frac{7x^{2} + 8xy - y^{2} + 23x - y - 18}{2x - 2y + 2}$ and $g(x,y)=8x + 2y - 4$. It is straightforward to verify that the inequality \(g(k,t) \leq n \leq 4k + 3t + 3\) implies \(k \leq \frac{t + 7}{4}\). Furthermore, we have \(g(k,t) \leq n \leq f(k, t)\). Therefore, based on the above definitions and by utilizing the properties of monotonic functions, together with straightforward calculations, it can be seen that the inequality \(\mu > \max\{\beta_1, \beta_2\}\) holds.
However, this leads to a contradiction between Equations~\eqref{seq2} and~\eqref{seq3}. This contradiction arises from our assumption that \( K_n \) does not contain a rainbow subgraph isomorphic to \( kP_4 \cup tP_2 \), while it does contain a rainbow subgraph \( H \cong (k-1)P_4 \cup (t+2)P_2 \).
Moreover, it is evident that for each $l \in [1, k-t+1]$, the equality $g(k,t)-g(k-l,t+2l) =2l $ holds. Additionally, it can be checked that if $k \geq t + 3$, then $f(k, t) \leq f(k-1, t+2)$.

Therefore, combining these observations with Theorem~\ref{th1}, we conclude that if $k \geq t$ and \(g(k,t) \leq n \leq f(k, t)\), then \( K_n \) contains a rainbow subgraph isomorphic to \( kP_4 \cup tP_2 \). This confirms Equation~\eqref{seq000} and completes the proof of the theorem. \qed

\section{Concluding remarks}
In this work, we explored the concept of neutral edges with respect to the anti-Ramsey number, a notion introduced by Erdős et al.\ in 1975. A subset of edges in a graph \(G\) is called \emph{neutral} if its removal does not affect the anti-Ramsey number of \(G\), which is the maximum number of colors in an edge-coloring of the complete graph $K_n$ that avoids a rainbow copy of \(G\).

Focusing on graphs of the form $G = kP_4 \cup tP_2$, where $S \subseteq E(G)$ consists of internal edges of the $P_4$ components, prior work \cite{GL-2} established that $S$ is neutral under any of the following conditions: (i) $t \geq k+1 \geq 2$ and $n \geq 8k + 2t - 4$; (ii) $k, t \geq 1$ and $n = 4k + 2t$; or (iii) $k = 1$, $t \geq 1$, and $n \geq 2t + 4$. However, neutrality does not hold when $k \geq 2$, $t = 0$, and $n = 4k$.

In this paper, we identify values of $k \geq t$ such that, for all $n$ in a specific subinterval of $[8k + 2t - 4, \infty)$, $S$ remains neutral. Since the anti-Ramsey numbers for matchings are well understood, this result provides a complete determination of the anti-Ramsey number for $G$ in these cases.

Based on our findings, we conjecture that this neutrality property may extend to the general case $t \geq 1$, $k \geq 1$, and $n \geq 4k + 2t$, but not to the case $t = 0$, $k \geq 2$, and $n \geq 4k$. Further investigation into the structure of neutral edge sets for other graph classes or parameter ranges remains an interesting direction for future research.
\vskip5mm
	
		
\noindent {\bf Acknowledgements} \vskip 3mm 

The research of Ali Ghalavand and Xueliang Li was supported the by the Natural Science Foundation of China under Grant No.12131013.  Jin was supported by the Natural Science Foundation of China under Grant No.12571380.
	
\vskip 3mm

\noindent {\bf Declaration of competing interest}
\vskip2mm

The authors declare that they have no known competing financial interests or personal relationships that could have appeared to influence the work reported in this paper.
\vskip2mm

	
	
\vskip 3mm

\noindent {\bf Data availability} 
	
No data was used in this investigation.


\begin{thebibliography}{99}

\bibitem{R2-SP-1} S. Akbari, A. Alipour, Multicolored trees in complete graphs, J. Graph Theory, 54 (2007), 221-232.

\bibitem{p1-1} A. Bialostocki, S. Gilboa, Y. Roditty, Anti-Ramsey numbers of small graphs, Ars Combin. 123 (2015), 41-53.

\bibitem{R2-SP-2} A. Bialostocki, W. Voxman, On the anti-Ramsey numbers for spanning trees, Bull. Inst.Combin. Appl. 32 (2001),  23-26.

\bibitem{p1-2} H. Chen, X.L. Li, J.H. Tu, Complete solution for the rainbow number of matchings, Discrete Math. 309 (10) (2009), 3370-3380.

 \bibitem{Er1} P. Erd\H{o}s, M. Simonovits, V.T. S\'{o}s, Anti-Ramsey theorems, Colloq. Math. Soc. J\'{a}nos Bolyai, Vol. 10 North-Holland Publishing Co., Amsterdam-London, 1975, pp. 633-643.

\bibitem{p1-4} C.Q. Fang, E. Gy\H{o}ri, M. Lu, J.M. Xiao, On the anti-Ramsey number of forests, Discrete Appl. Math. 291 (11) (2021), 129-142.

\bibitem{p1-5} S. Fujita, A. Kaneko, I. Schiermeyer, K. Suzuki, A rainbow $k$-matching in the complete graph with $r$ colors, Electron. J. Combin. 16 (2009), \#R51.

\bibitem{GL-2} A. Ghalavand, Q. Jie, Z. Jin, X. Li, L. Pan, On the anti-Ramsey number under edge deletion, arXiv: 2511.06034 [math.CO] (8 Nov 2025).

\bibitem{GL-1} A. Ghalavand, X. Li, On the anti-Ramsey number of spanning linear forests with paths of lengths 2 and 3, arXiv: 2509.25949 [math.CO] (30 Sep 2025).

\bibitem{p1-6} S. Gilboa, Y. Roditty, Anti-Ramsey numbers of graphs with small connected components, Graphs Combin. 32 (2) (2016), 649-662.

\bibitem{R2-HP-1} R. Gu, J. Li, Y. Shi, Anti-Ramsey numbers of paths and cycles in hypergraphs, SIAM J. Discrete Math. 34 (1) (2020), 271-307.

\bibitem{R2-HP-2} M. Guo, H. Lu, X. Peng, Anti-Ramsey number of matchings in 3-uniform hypergraphs, SIAM J. Discrete Math. 37 (3) (2023), 1970-1987.

\bibitem{p1-7} R. Haas, M. Young, The anti-Ramsey number of perfect matching, Discrete Math. 312 (2012), 933-937.

\bibitem{p1-8} M.L. He, Z.M. Jin, Rainbow short linear forests in edge-colored complete graph, Discrete Appl. Math. 361 (2025), 523-536.

\bibitem{p1-9} S. Jahanbekam, D.B. West, Anti-Ramsey problems for t edge-disjoint rainbow spanning subgraphs: cycles, matchings, or trees, J. Graph Theory 82 (1) (2016), 75-89.

\bibitem{new-3} Q. Jie, Z.M. Jin, Anti-Ramsey number of union of 5-path and matching, Discuss.  Math. Graph Theory 45 (2025), 1185-1210.

\bibitem{p1-10} Q. Jie, M.L. He, Z.M. Jin, Rainbow forest consisting of short paths in $K_n$, Discrete Appl. Math. 376 (2025), 260-269.

\bibitem{new-2} Q. Jie, Z.M. Jin, Rainbow-free colorings for spanning linear forest consisting of short paths, 2025, {\it submitted}.

\bibitem{new-1} Z.M. Jin and J.Q. Gu, Rainbow disjoint union of clique and matching in edge-colored complete graph, Discuss. Math. Graph Theory 44 (2024), 953–970.

\bibitem{n-p4-1} Z. Jin, Q. Jie, Z. Cao, Rainbow disjoint union of $P_4$ and a matching in complete graphs, Appl. Math. Comput. 474 (2024), 128679.

\bibitem{R2-SP-3} L.Y. Lu, A. Meier, Z.Y. Wang, Anti-Ramsey number of edge-disjoint rainbow spanning trees in all graphs, SIAM J. Discrete Math. 37 (2) (2023), 1162-1172.

\bibitem{R2-SP-4} L.Y. Lu, Z.Y. Wang, Anti-Ramsey number of edge-disjoint rainbow spanning trees, SIAM J. Discrete Math. 34 (4) (2020), 2346--2362.

\bibitem{p3-6} J.J. Montellano-Ballesteros, V. Neumann-Lara, An anti-Ramsey theorem on cycles, Graphs Combin. 21 (3) (2005), 343-354.

 \bibitem{p2-9} I. Schiermeyer, Rainbow numbers for matchings and complete graphs, Discrete Math. 286 (2004), 157-162.

\bibitem{p2-10} M. Simonovits, V.T. Sós, On restricted coloring of $K_n$, Combinatorica 4 (1) (1984), 101–110.

\bibitem{R2-HP-3} Y. Tang, T. Li, G. Yan, Anti-Ramsey number of expansions of paths and cycles in uniform hypergraphs, J. Graph Theory 101 (4) (2022), 668-685. 

\bibitem{p1-14} T.Y. Xie, L.T. Yuan, On the anti-Ramsey numbers of linear forests, Discrete Math. 343(12)(2020), 112130.


\end{thebibliography}
\end{document}